\documentclass[11pt, reqno]{amsart}
\usepackage[english]{babel}
\usepackage{amsmath,amssymb,graphicx,enumerate, 
xcolor
}
\usepackage{amstext}
\usepackage{amsthm}
\usepackage{MnSymbol}
\usepackage{wasysym}

\usepackage[
hypertexnames=false, colorlinks, citecolor=black, linkcolor=blue, urlcolor=red]{hyperref}

\newcommand{\mathd}{\mathrm{d}}
\newcommand{\mathi}{\mathrm{i}}
\newcommand{\mathc}{\mathrm{c}}

\newcommand{\tmop}[1]{\ensuremath{\operatorname{#1}}}

\newcommand{\bI}{\mathbf{I}}
\newcommand{\1}{\mathbf{1}}
\newcommand{\bM}{\mathbf{M}}

\newcommand{\R}{\mathbb{R}}

\newcommand{\cD}{\mathcal{D}}

\newcommand{\cS}{\mathcal{S}}
\newcommand{\cH}{\mathcal{H}}

\newcommand{\lla}{\llangle}
\newcommand{\rra}{\rrangle}

\newcommand{\lnm}{\left|}
\newcommand{\rnm}{\right|}

\newcommand{\e}{\varepsilon}

\newcommand{\f}{\varphi}
\newcommand{\fdot}{\,\cdot\,}

\newcommand{\wt}{\widetilde}
\newcommand{\tr}{\operatorname{trace}}

\numberwithin{equation}{section}

\newtheorem{theorem}{Theorem}[section]

\newtheorem{lemma}[theorem]{Lemma}

\newtheorem*{prop*}{Proposition}

\theoremstyle{remark}
\newtheorem{remark}[theorem]{Remark}
\newtheorem*{rem*}{Remark}

\newcommand{\ci}[1]{_{ {}_{\scriptstyle #1}}}

\newcommand{\ti}[1]{_{\scriptstyle \text{\rm #1}}}
\newcommand{\ut}[1]{^{\scriptstyle \text{\rm #1}}}

\newenvironment{entry}
{\begin{list}{X}%
		{%
			\setlength{\labelwidth}{55pt}%
			\setlength{\leftmargin}{\labelwidth}
			\addtolength{\leftmargin}{\labelsep}%
			\setlength{\itemsep}{.4pc}%
		}%
	}%
	{\end{list}}   

\renewcommand{\labelenumi}{\textup{(\roman{enumi})}}

\newcounter{vremennyj}

\newcommand\cond[1]{\setcounter{vremennyj}{\theenumi}\setcounter{enumi}{#1}\labelenumi\setcounter{enumi}{\thevremennyj}}

\begin{document}

\title{The matrix-weighted dyadic convex body  maximal operator  is not bounded}

\author{F. Nazarov, S. Petermichl, K. A. \v{S}kreb, S. Treil}

\begin{abstract}
The convex body maximal operator is a natural generalization of the Hardy--Littlewood maximal operator. In this paper we are considering its dyadic version in the presence of a matrix weight.  To our surprise it turns out that this operator is not bounded. This is in a sharp contrast to a Doob's inequality in this context.  At first, we show that the convex body Carleson Embedding Theorem with matrix weight fails. We then 
  deduce the unboundedness of the matrix-weighted convex body maximal operator. 
\end{abstract}

\maketitle

\setcounter{section}{-1}
\section{Notation}

\begin{entry}
\item[$I_0$] interval $[0, 1]$;

\item[$\mathcal{D}$] the dyadic lattice, i.e.\hspace{-0.05em}, the collection of all dyadic intervals; 

\item[$I_+$, $I_-$] left and right halves of the interval $I\in\mathcal{D}$;

\item[$\mathcal{D}_{\pm}$] $\mathcal{D}_{\pm}=\{I_\pm : I \in \mathcal{D}\}$ so that $\mathcal{D}=\mathcal{D}_+ \dot{\cup} \mathcal{D}_- \dot{\cup} \{I_0\}$;

\item[$\mathcal{D}(K)$] $\mathcal{D} (K) = \{ I\in\mathcal{D} : I\subset K\}$ for $K\in\mathcal{D} $;

\item[$| I |$] the Lebesgue measure of the set $I\subset I_0$; 

\item[$\mathcal{D}^n$] $\mathcal{D}^n:= \{I\in\mathcal{D}: |I| = 2^{-n}\}$;

\item[$\mathcal{D}^{\leqslant n}$] $\mathcal{D}^{\leqslant n} := \bigcup_{k\leqslant n} 
\mathcal{D}^k 
=\{I\in\mathcal{D}: |I|\geqslant 2^{-n} \}$;

\item[$\mathcal{C}_\pm$] for a collection 
$\mathcal{C}\subset \mathcal{D}$ of dyadic intervals, $\mathcal{C}_\pm = \mathcal{C} \cap \mathcal{D}_\pm$;

\item[$(\mathcal{F}_n)_{n\geqslant 0}$]  the dyadic filtration: $\mathcal{F}_n$ is the 
$\sigma$-algebra 
generated by $\mathcal{D}^n$; 

\item[ $\langle f \rangle\ci I$] average of a scalar, vector or matrix function $f$ 
over $I$: $\langle f \rangle\ci I := | I |^{- 1} \int_I f (x) \mathd x$; 

\item[$f (I)$] for $I\subset I_0$ we denote by $f(I)$ the integral of a 
scalar, vector or matrix function $f$ over $I$: $f
(I) = \int_I f (x) \mathd x$; 

\item[$\langle \fdot, \fdot \rangle\ci {\mathbb{R}^d}$] the standard inner product in
$\mathbb{R}^d$; 

\item[$\lnm \fdot \rnm\ci{\mathbb{R}^d}$]  the norm induced by $\langle \cdot, \cdot
\rangle\ci {\mathbb{R}^d}$ with subscript omitted if $d=1$; 

\item[$\lnm\fdot\rnm\ti{op}$] the operator norm of a matrix; 

\item[$\| \fdot \|\ci{X}$] norm in the function space $X$; 

\item[$B(X)$] unit ball in the normed space $X$; 


\item[ $\1\ci I$] the characteristic function of $I$.
\end{entry}

\

Generally, since we are dealing with vector- and matrix-valued functions we
will use the symbol $\| \fdot \|$ (usually with a subscript) for the norm in a
functions space, while $\lnm \fdot \rnm$ is used for the norm in the underlying
vector (matrix) space. Thus, for a vector valued function $f$, the symbol
$\| f \|_{L^2}$ denotes its $L^2$-norm, but the symbol $\lnm f \rnm$ stands for the
scalar valued function $x \mapsto \lnm f (x) \rnm$. 

Finally, we will use the \emph{linear algebra notation}, identifying vector $a$ in a Hilbert space 
$\cH$ with the operator  $\alpha\mapsto \alpha a$ acting from scalars to $\cH$. In this case the 
symbol $a^*$ denotes the (bounded) linear functional $x\mapsto \langle x, a \rangle $. In our case $\cH$ is the 
real Hilbert space $\R^d$, vectors in $\R^d$ are the column vectors, and $a^*$ is just the transpose 
of the column $a$.

\section{Introduction}

The simple dyadic maximal function
\[ \mathcal{M}f (x) = \sup_{I \in \mathcal{D}: x \in I} \lnm \langle f \rangle\ci I \rnm \]
and the dyadic Hardy-Littlewood maximal function
\[ \mathcal{M}^{\mathc} f (x) = \sup_{I \in \mathcal{D}: x \in I} 
\langle \lnm f \rnm \rangle\ci I \]
together with their scalar weighted analogues
\[ \mathcal{M}_w f (x) = \sup_{I \in \mathcal{D}: x \in I} \frac{1}{w (I)}
   \lnm \int_I f (y) w(y) \mathd y \rnm \]
and
\[ \mathcal{M}^{\mathc}_w f (x) = \sup_{I \in \mathcal{D}: x \in I} \frac{1}{w (I)}
   \int_I \lnm f (y) \rnm w(y) \mathd y \]
for positive $w \in L^1$ are classical objects in harmonic analysis. (The reason for the superscript $\mathc$ will become clear in the next section.)

The inequality 
$$\| \mathcal{M}_w \|\ci {L_w^2 (\R) \rightarrow L_w^2(\R)} \leqslant 2$$ 
is a special case of Doob's inequality and the
same inequality for $\mathcal{M}^{\mathc}_w$ follows immediately from the
observations that $\mathcal{M}_w^{\mathc} f =\mathcal{M}_w \lnm f \rnm$ and $\| \lnm f \rnm
\|\ci {L_w^2 (\R)} = \| f \|\ci {L_w^2 (\R)}$.

The boundedness of $\mathcal{M}_w$ or $\mathcal{M}^{\mathc}_w$ in $L_w^2
(\R)$ can be rewritten as the boundedness of the operators
\[ M_w f (x) = \sup_{I \in \mathcal{D}: x \in I} \left| w^{1 / 2} (x)
   \frac{1}{w (I)} \int_I f (y) w (y) \mathd y \right| {\color{red} } \]
and
\[ M^{\mathc}_w f (x) = \sup_{I \in \mathcal{D}: x \in I} w^{1 / 2} (x) \frac{1}{w
   (I)} \int_I | f (y) | w (y) \mathd y \]
from $L_w^2 (\R)$ to the usual $L^2 (\R)$.

The natural counterpart (see the next section for details) of $M_w f$ and $M^{\mathc}_w f$ in the case when $f$ is a vector-valued
function and $W$ a matrix weight (a positive definite matrix function) are
\[ M\ci W f (x) = \sup_{I \in \mathcal{D}: x \in I} 
\lnm W^{1 / 2} (x) \langle W
   \rangle^{- 1}\ci I \langle W f \rangle\ci I \rnm\ci {\R^d} \]
and
\[ M^{\mathc}\ci W f (x) = 
\sup_{\substack{I\in\cD\,:\, x\in I \\ \f_I : I\to [-1, 1]}} 
\lnm W^{1 / 2} (x) \langle W \rangle^{- 1} \ci I \langle \varphi\ci I W f
   \rangle \ci I \rnm \ci {\R^d} . \]
When $d=1$ use the signum of $Wf$ on $I$ for $\varphi\ci I$  to get the usual absolute value.  
This definition will be shown in the next section to be equivalent to a Christ-Goldberg type definition.   
The boundedness of these operators from $L\ci W^2
(\R^d)$ to $L^2 (\R)$ when $W = \tmop{Id}_d$ remains a
simple consequence of Doob's inequaliy, but is more complicated for general
matrix weights. The positive
result for $M\ci W f$ was established in \cite{PePoRe18} by reducing it to the weighted
Carleson Embedding Theorem from \cite{CuTr15}. :

\begin{theorem}\label{WCET}
  (Culiuc, Treil) Let $W$ be a matrix weight and $(A\ci I)\ci {I \in \mathcal{D}}$ a
  sequence of positive definite matrices. Then the
  following are equivalent:
  \begin{enumerate}
    \item There exists $c_{\cond1} > 0$ such that for all $K \in \mathcal{D}$,
    \[ \frac{1}{| K |} \sum_{I \in \mathcal{D} (K)} \langle W \rangle\ci I A\ci I
       \langle W \rangle\ci I \leqslant c_{\cond1} \langle W \rangle\ci K . \]
    \item  There exists $c_{\cond2} > 0$ such that for all $f \in L\ci W^2
    (\R^d)$,
    \[ \sum_{I \in \mathcal{D}} \lnm A^{1 / 2}\ci I \langle W f \rangle\ci I
       \rnm\ci {\R^d}^2 \leqslant c_{\cond2} \| f \|\ci {L\ci W^2 (\R^d)}^2
       . \]
  \end{enumerate}
  Moreover, the best possible constants $c_{(\mathi)}$ and $c_{(\mathi \mathi)}$ in $(\mathi)$ and $(\mathi \mathi)$
  satisfy $c_{(\mathi)} \leqslant c_{(\mathi \mathi)} \leqslant C c_{(\mathi)}$ with $C > 0$ only depending
  on the dimension $d$ but not on $W$. 
\end{theorem}

The bound for the norm $\| M\ci W \|\ci {L\ci W^2 (\R^d) \rightarrow L^2
(\R)}$ follows by the so-called linearization technique and the exact
statement is as follows:

\begin{theorem}
  (Petermichl, Pott, Reguera) Let $W$ be a matrix weight. Then
  \[ \| M\ci W \|\ci {L\ci W^2 (\R^d) \rightarrow L^2 (\R)} \leqslant
     C, \]
  where $C$ depends only on $d$. 
\end{theorem}

These results raised hopes that the larger maximal
function $M^{\mathc}\ci W$ may be bounded from $L\ci W^2 (\R^d)$ to $L^2
(\R)$ as well, so it came as a surprise to us that in general it
actually is not. The main purpose of this paper is to present and discuss the
counter-example. Again, to construct it, we show
that the corresponding convex body Carleson Embedding Theorem fails, which was
also quite surprising for us and which may be of independent interest. Here are 
the precise formulations of our negative results:

\begin{theorem} \label{noWcCET}
  For any fixed dimension $d \geqslant 2$, there exist a matrix weight $W$ and
  a sequence $(A\ci I)\ci {I \in \mathcal{D}}$ of positive definite matrices such that there exists $c_{\cond1}> 0$ such that for all $K \in
  \mathcal{D} $
  \begin{equation} \label{e: Carl test 02}
  \frac{1}{| K |} \sum_{I \in \mathcal{D} (K)} \langle W \rangle\ci I A\ci I
     \langle W \rangle\ci I \leqslant c_{\cond1} \langle W \rangle\ci {K}, 
 \end{equation}
  but there exist $f \in L\ci W^2 (\R^d)$ and a sequence
  $(\varphi\ci I)\ci {I \in \mathcal{D}} $, $- 1 \leqslant
  \varphi\ci I \leqslant 1$, such that
  \begin{equation}\label{e: Carl blow-up} 
  \sum_{I \in \mathcal{D} } \lnm A^{1 / 2}\ci I \langle \varphi\ci I W f
     \rangle\ci I \rnm\ci {\R^d}^2 = \infty . 
  \end{equation}
\end{theorem}

\begin{theorem} \label{noWcMAX}
  For any fixed dimension $d \geqslant 2$, there exists a matrix weight $W$ such
  that $M^{\mathc}\ci W$ does not map $L\ci W^2 (\R^d)$ to $L^2 (\R)$. 
\end{theorem}

\section{Discussion of definitions}

Let $f$ be a vector-valued function with values in $\R^d$. A
Hardy-Littlewood maximal function of $f$ at $x$ is a quantity
allowing one to control all averages of $f$ over
intervals containing $x$. The most straightforward quantity of that type is
just $\sup_{I \in \mathcal{D}: x \in I} \langle \lnm f \rnm\ci {\R^d}
\rangle\ci I$. However, this quantity is often too crude when $f$ is large in
some directions and small in some other ones (all such information is
completely lost here) and is too strongly tied to the particular choice of the
Euclidean norm in $\R^d$ to yield useful bounds in
the case when the size of $f$ is measured in some other norm, especially in a
norm depending on a point, as it is the case in
matrix-weighted spaces.

A more reasonable idea would be to use the convex-body average $\lla
 f \rra \ci I$ of an $L^1$ function $f$ over an interval $I$, defined as
\[ \lla f \rra\ci I = \{ \langle \varphi\ci I f \rangle\ci I :\;\;
   \varphi\ci I : I \rightarrow [- 1, 1] \} . \]
$\lla f \rra \ci I$ is always a symmetric compact convex set
in $\R^d$ (the compactness follows from the weak-$\ast$ compactness
of the unit ball in $L^{\infty} (I,\R) = L^1 (I,\R)^{\ast}$, where these spaces consist of functions defined on the interval $I$) but, in general, it may contain no inner points.
However, contrary to the terminology of some convex body geometry books, we
will still call it a convex body.

This convex body average makes sense even if $d = 1$, in which case it is
just the interval $[- \langle \lnm f \rnm \rangle\ci I, \langle \lnm f \rnm \rangle\ci I]$, so there it carries just as
much information as $\langle \lnm f \rnm \rangle\ci I$. \

Let now $\rho$ be any norm in $\R^d$. For a set $\Omega \subset
\R^d$, define its norm as
\[ \rho (\Omega) = \sup_{x \in \Omega} \rho (x) . \]
If $\rho (x) = \lnm x \rnm\ci {\R^d}$ is the usual Euclidean norm, then we
have 
\[ d^{- 1} \langle \lnm f \rnm \ci {\R^d} \rangle\ci I \leqslant \rho (\lla
    f \rra \ci I) \leqslant \langle \lnm f \rnm\ci {\R^d}
   \rangle\ci I . \]
The right inequality is just the triangle inequality. To obtain the left one,
just write $f = (f_1, \ldots, f_d)$ and note that for all $ k$
\[ \langle \lnm f \rnm\ci {\R^d} \rangle\ci I \leqslant \sum^d_{k = 1} \langle \lnm
   f_k \rnm \rangle\ci I, \]
so there exists $k$ with $\langle \lnm f_k \rnm \rangle\ci I \geqslant d^{- 1} \langle
\lnm f \rnm\ci {\R^d} \rangle\ci I$. On the other hand, choosing $\varphi\ci I =
\tmop{sign} f_k$, we get
\[ \rho (\lla f \rra \ci I) \geqslant (\langle \varphi\ci I f 
   \rangle\ci I)_k = \langle \lnm f_k \rnm \rangle\ci I \]
and the left inequality follows.

The advantage of the convex body approach is that this inequality is preserved
if we change the standard Euclidean norm in $\R^d$ to any other
Euclidean norm, i.e., if we take a positive definite matrix $A$ and consider
$\rho_A (x) = \lnm A^{1 / 2} x \rnm\ci {\R^d}$. The corresponding estimates
\begin{equation}\label{sim}
 d^{- 1} \langle \rho\ci A (f) \rangle\ci I \leqslant \rho\ci A (\lla f
   \rra\ci I) \leqslant \langle \rho\ci A (f) \rangle\ci I \end{equation}
immediately follow from the standard Euclidean ones by considering $A^{1 / 2}
f$ instead of $f$.

Inspired by the idea of the convex body average, we can now define a convex
body valued maximal function.

The simple maximal function $\bM$ will be just
\[ \bM f (x) = \left\{ \sum_{I \in \mathcal{D}: x \in I} a\ci I \langle f
   \rangle\ci I :\;\; a\ci I \in \R, \sum_{I \in \mathcal{D}: x \in I} 
   | a\ci I | \leqslant 1 \right\} \]
i.e., $\bM f (x)$ is the absolute convex hull of averages of
$f$ over intervals containing $x$. The generalization of the Hardy-Littlewood
maximal function will now be
\[ \bM^{\mathc} f (x) = \left\{ \sum_{I \in \mathcal{D}: x \in I} a\ci I \lla
   f \rra\ci I :\;\; a\ci I \in \R, \sum_{I \in \mathcal{D}:
   x \in I} \lnm a\ci I \rnm \leqslant 1 \right\}, \]
where the sum of convex bodies is understood in the Minkowski sense: $A + B =
\{ a + b : \;\; a \in A, b \in B \}$. Plugging in the definition of 
$\lla f \rra \ci I$, we can also rewrite this as
\[ \bM^{\mathc} f (x) = \left\{ \sum_{I \in \mathcal{D}: x \in I} a\ci I \langle
   \varphi \ci I f \rangle\ci I :\;\; a\ci I \in \R, \sum_{I \in
   \mathcal{D}: x \in I} | a\ci I | \leqslant 1, \; \varphi\ci I : I \rightarrow [- 1,
   1] \right\} . \]
For a matrix weight $W${\hspace{-0.2em}}, one can easily write the
weighted analogues $\bM\ci W$ and $\bM\ci W^{\mathc}$ of $\bM$ and $\bM^{\mathc}$:
\[ \bM\ci W f (x) = \left\{ \sum_{I \in \mathcal{D}: x \in I} a\ci I \langle W
   \rangle^{- 1}\ci I \langle W f \rangle\ci I :\;\; a\ci I \in \R, \sum_{I
   \in \mathcal{D}: x \in I} | a\ci I | \leqslant 1 \right\} \]
and
\[ \bM\ci W^{\mathc} f (x) = \left\{ \sum_{I \in \mathcal{D}: x \in I} a\ci I \langle W
   \rangle^{- 1}\ci I \langle \varphi\ci I W f \rangle\ci I :\;\; a\ci I \in
   \R, \sum_{I \in \mathcal{D}: x \in I} | a\ci I | \leqslant 1, \; 
   \varphi\ci I : I \rightarrow [- 1, 1] \right\} . \]
The values here are again convex bodies.
 
The weighted space $L\ci W^2 (\R^d)$ is just the space of all measurable
functions $f$ with values in $\R^d$ for which
\[ \| f \|_{L\ci W^2 (\R^d)}^2 = \int \lnm W^{1 / 2} f \rnm\ci {\R^d}^2
   = \int \langle W f, f \rangle\ci {\R^d} < + \infty . \]
In other words, it is the $L^2$ space in which the size of $f$ is measured not
by the usual Euclidean norm in $\R^d$ but by the norm $\rho\ci {W (x)}$.
In line with our definition at the beginning of this section, we can now write
\[ \| \bM\ci W f \|_{L\ci W^2 (\R^d)}^2 = \int \rho\ci {W (x)} (\bM\ci W f
   (x))^2 \mathd x \]
and similarly for $\bM^{\mathc}\ci W$, where, as before, the norm of a subset of
$\R^d$ is understood as the supremum of the norms of the vectors
contained in it. Since $\bM\ci W f (x)$ is the absolute convex hull of $\langle W
\rangle^{- 1}\ci I \langle W f \rangle\ci I$ with $I \in \mathcal{D}$ such that $x \in I$, by
the convexity property of the norm,
\[ \rho\ci {W (x)} (\bM\ci W f (x)) = \sup_{I \in \mathcal{D}: x \in I} \rho\ci {W (x)}
   (\langle W \rangle^{- 1}\ci I \langle W f \rangle\ci I) = \sup_{I \in
   \mathcal{D}: x \in I} \lnm W^{1 / 2} (x) \langle W \rangle^{- 1}\ci I \langle W f
   \rangle\ci I \rnm\ci {\R^d} \]
and, by the same logic,
\[ \rho\ci {W (x)} (\bM^{\mathc}\ci W f (x)) =\sup_{\substack{I\in\cD\,:\, x\in I \\ \f_I : I\to [-1, 1]}}  
 \lnm W^{1 / 2} (x) \langle W \rangle^{- 1}\ci I \langle
   \varphi\ci I W f \rangle\ci I \rnm\ci {\R^d} \]
leading to the maximal functions $M\ci W$ and $M\ci W^{\mathc}$ used in the introduction so
that
\[ \| \bM\ci W f \|\ci{L\ci W^2 (\R^d)} = \| M\ci W f \|\ci{L^2 (\R)} \]
and similarly for $\bM^{\mathc}\ci W$ and $M^{\mathc}\ci W$.

The reader familiar with Christ-Goldberg definitions of maximal functions in
the context of matrix weighted spaces may prefer to consider the quantity
\[ \frac{1}{| I |} \int_I \lnm W^{1 / 2} (x) \langle W \rangle^{- 1}\ci I W (y) f
   (y) \rnm\ci{\R^d} \mathd y \]
instead of our
\[ \sup_{\varphi\ci I : I \rightarrow [- 1, 1]} \lnm W^{1 / 2} (x) \langle W
   \rangle^{- 1}\ci I \langle \varphi_I W f \rangle\ci I \rnm\ci{\R^d} . \]
However, our inequalities in (\ref{sim}) show that these
quantities are the same up to a factor of $d$.

The final aspect we want to discuss is the
measurability of $\bM\ci W$ and $\bM^{\mathc}\ci W$ and another natural way to define the
$L\ci W^2 (\R^d)$ norm of set-valued functions. We shall do this
discussion for $\bM\ci W^{\mathc}$. The case of $\bM\ci W$ is the same but simpler.

Let us introduce the truncated maximal functions $M^{\mathc}\ci {W, n}$ and $\bM^{\mathc}\ci {W,
n}$ in which we take into consideration only the intervals $I \in
\mathcal{D}^{\leqslant n}$, so
\[ \bM\ci {W, n}^{\mathc} f (x) = \left\{ \sum_{I \in \mathcal{D}^{\leqslant n}:x \in I} 
a\ci I
   \langle W \rangle^{- 1}\ci I \langle \varphi\ci I W f \rangle\ci I :\;\; a\ci I \in
   \R, \sum_{I \in \mathcal{D}^{\leqslant n} : x \in I} | a\ci I |
   \leqslant 1,  \; \varphi\ci I : I \rightarrow [- 1, 1] \right\} \]
and 
\[ M^{\mathc}\ci{W, n} f (x) = \sup_{v \in \bM_{W, n}^{\mathc} f (x)} 
\lnm W^{1 / 2} (x) v \rnm\ci{\R^d} = 
   \sup_{\substack{I\in\cD\,:\, x\in I \\ \f_I : I\to [-1, 1]}} 
   \lnm W^{1 / 2} (x) \langle W \rangle^{-
   1}\ci I \langle \varphi\ci I W f \rangle\ci I \rnm\ci{\R^d} . \]
Notice that $\bM^{\mathc}\ci{W, n} f (x)$ is a symmetric compact convex set. Moreover,
this set depends only on the $n$-th level dyadic interval $I \in
\mathcal{D}^n$ the point $x$ belongs to.

Let now $K \subset \R^d$ be any convex compact set and let $\varepsilon >
0.$ Then there exists a finite set $\{ v_1, \ldots ,v_N \} \subset K$ such that
$K \subset (1 + \varepsilon) \tmop{conv} (v_1, \ldots ,v_N)$. In particular, it
implies that for every $x \in I,$ 
\[ \max_{1 \leqslant j \leqslant N} \lnm W^{1 / 2} (x) v_j \rnm\ci {{\R}^d}
   \leqslant \sup_{v \in K} \lnm W^{1 / 2} (x) v \rnm\ci{\R^d} 
   \leqslant (1 +\varepsilon) \max_{1 \leqslant j \leqslant N} \lnm W^{1 / 2} (x) v_j
   \rnm\ci{\R^d} . \]
For each $x \in I$, define $j (x)$ as the least index $j$ for which $\lnm W^{1 /
2} (x) v_j \rnm\ci{\R^d}$ is maximal. In other words, $j (x) = J$ if and
only if $\lnm W^{1 / 2} (x) v_j \rnm\ci{\R^d} < \lnm W^{1 / 2} (x) v\ci J
\rnm\ci{\R^d}$ for $j < J$ and $\lnm W^{1 / 2} (x) v_j \rnm\ci{\R^d}
\leqslant \lnm W^{1 / 2} (x) v\ci J \rnm\ci{\R^d}$ for $j > J$. Then $g (x) =
v_{j (x)}$ is a measurable function on $I$ such that $g (x) \in K$ and
\[ \lnm W^{1 / 2} (x) g (x) \rnm\ci{\R^d} \geqslant (1 + \varepsilon)^{- 1}
   \sup_{v \in K} \lnm W^{1 / 2} (x) v \rnm\ci{\R^d} \]
for all $x \in I$. 

Applying this construction to every $I \in \mathcal{D}^n$ with $K = \bM\ci{W,
n}^{\mathc} f |\ci I$, we get a measurable $g : I_0 \rightarrow \R^d$
such that $g (x) \in \bM\ci {W, n}^{\mathc} f (x)$ and
\[ \lnm W^{1 / 2} (x) g (x) \rnm\ci{\R^d} > (1 + \varepsilon)^{- 1} M^{\mathc}\ci{W,
   n} f (x) \]
for all $x \in I_0$. Letting $\varepsilon \rightarrow 0$ (along some
sequence), we obtain that $M^{\mathc}\ci {W, n} f$ is the supremum of a countable family
of measurable functions and, thereby, measurable. Finally, letting $n
\rightarrow \infty$, we get that $M^{\mathc}\ci W f$ is measurable and, moreover,
\[ \| M^{\mathc}\ci W f \|\ci{L^2 (\R)} = \sup \{ \| g \|\ci {L\ci W^2 (\R^d)}
   : \;\; g : I_0 \rightarrow \R^d, \; g \, \tmop{ is } \, \tmop{ measurable}
   , \; g \in \bM \ci W^{\mathc} f \} . \]

\

\section{Proofs}

We first prove Theorem \ref{noWcCET}, which states that the convex body analogue of the weighted 
Carleson Embedding Theorem (i.e.\hspace{-0.05 em}, of Theorem \ref{WCET}) fails. 

Using this result we then conclude the failure of the maximal function estimate.

\subsection{Failure of the convex body Carleson Embedding Theorem: proof of Theorem \ref{noWcCET}}
\label{s: fail CBCET}

\subsubsection{Overview of the construction}
\label{s: overview no CET}

We consider a counterexample in $\R^2$, which trivially extends to an example in all higher dimensions.

We will be constructing the weight $W$ supported on the interval $I_0$ as a martingale, from the top down. This means at the  $n$th step 
($n=0, 1, 2, \ldots$)  we 
construct the averages $W\ci I$, $I\in\cD^n$, which will be the averages of 
the constructed weight $W$. The averages $W\ci I =\langle W \rangle \ci I$ should satisfy the \emph{martingale dynamics}, 
\begin{align}
\label{e: mart dyn}
W\ci I = \frac12\left( W\ci{I_+} + W\ci{I_-}\right) , 
\end{align}
so the sequence of weights $W_n$, 
\begin{align*}
W_n := \sum_{I\in\cD^n} W\ci I \1\ci I
\end{align*}
is a martingale with respect to the dyadic filtration. In our case the sequence $W_n$ will be 
uniformly bounded, so we will have convergence%
\footnote{For the convergence a much weaker condition on \emph{uniform integrability} is 
sufficient} 
$W_n\to W$ (say, entry-wise in $L^1$), 
\begin{align*}
W\ci I = \langle W \rangle\ci I\qquad \forall I \in\cD. 
\end{align*}

We will work with the spectral decompositions of the averages $W\ci I$
\begin{align*}
W\ci I = \alpha\ci I a\ci I a\ci I^{\ast} + \beta\ci I b\ci I b\ci I^{\ast}
\end{align*}
where $\alpha\ci I$, $\beta\ci I$ are the eigenvalues, and $a\ci I$, $b\ci I$ are the corresponding 
normalized eigenvectors of $W\ci I$. This means, in particular, that $| a\ci I |\ci{\R^2} = |
b\ci I | \ci{\R^2} =1$ and $a\ci I \perp b\ci I$.  

In our situation, the eigenvalues $\alpha\ci I$, $\beta\ci I$ will depend only on $|I|$, i.e.\hspace{-0.05 em}, for 
$I\in\cD^n$ we will have $\alpha\ci I = \alpha_n$, $\beta\ci I =\beta_n$. Note that this condition 
means  that $\tr W\ci I =\tr W\ci J$ if $I, J\in \cD^n$; then the martingale property \eqref{e: 
mart dyn} and linearity of the trace imply that $\tr W\ci I =\tr W\ci{I_0}$ for all $I\in\cD$, 
which, in turn, implies the uniform boundedness of $W_n$, $\lnm W_n (x) \rnm\ti{op} \leqslant 
\alpha_0+\beta_0$. 

In our construction, we will have the condition numbers (eccentricity) $\alpha_n/\beta_n \to\infty$  
as $n\to \infty$. The operators $A\ci I$ will be just appropriately renormalized projections onto 
$b\ci I$, so the largest term $\alpha\ci I a\ci I a\ci I^*$ of $W\ci I$ disappears from the testing 
condition (\ref{e: Carl test 02}) in  Theorem \ref{noWcCET}. 

However,  an appropriate (and very simple) choice of the functions $\f\ci I$ in the conclusion  (\ref{e: Carl blow-up}) of  Theorem \ref{noWcCET} 
 will bring this huge term $\alpha\ci I a\ci I a\ci I^*$ into play, and that will give us the desired blow-up.

\begin{figure}[h]\centering
\includegraphics{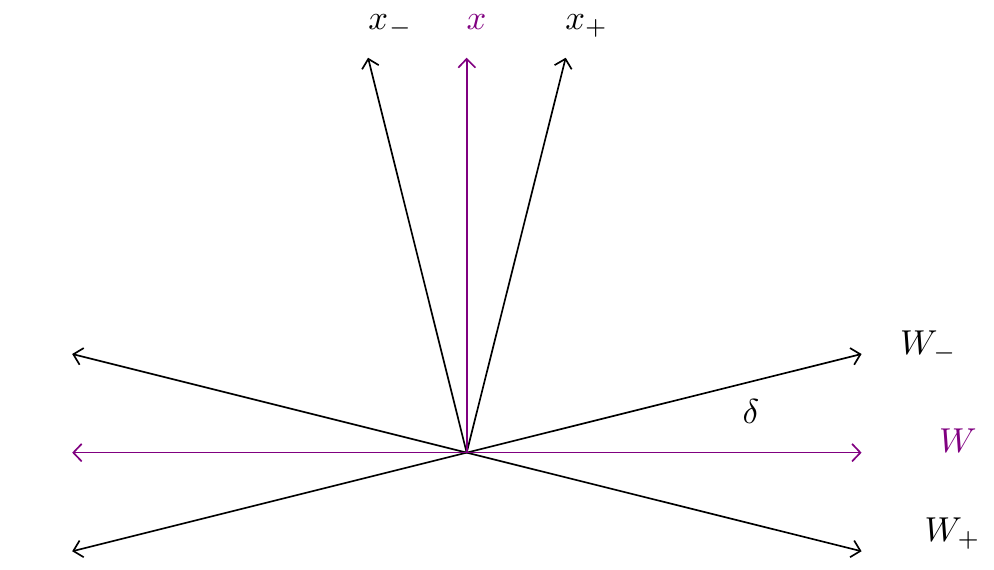}
\caption{In this picture we see the images of the unit ball under the weight $W$,
		sketched as arrows though they are thin ellipses. Further, we see the
		directions of the vectors $b$.}
	\label{fig: W_I}
\end{figure}



\subsubsection{Gory details}

Now let us assume that for all $I\in\cD^n$, we constructed $W\ci I$, and its spectral decomposition 
is given by 
\[  
W \ci I = \alpha_n a\ci I a^{\ast}\ci I + \beta_n b\ci I b^{\ast}\ci I  .
\]
Recall that $a\ci I$ and $b\ci I$ is an orthonormal pair of vectors. Let us define the 
averages $W\ci{I_\pm}$ as 
\[ 
W \ci{I_{\pm}} = \alpha_{n + 1} a\ci{I_{\pm}}  a\ci{I_{\pm}}^{\ast} + \beta_{n + 1} b\ci{I_{\pm}} 
b^{\ast}\ci{I_{\pm}}, 
\]
where $a\ci{I_{\pm}}$ and  $b\ci{I_{\pm}}$ are small rotations of the vectors $a\ci I$ and $b\ci 
I$, namely, for some small $\delta_{n+1}>0$ 
%
\begin{align}
\label{e: ab to+-}
a\ci{I_{\pm}} = \left( \frac{1}{1 + \delta_{n + 1}^2} \right)^{1/2}
(a\ci I \pm \delta_{n + 1} b\ci I),  \qquad b\ci{I_{\pm}} = \left( \frac{1}{1 +
	\delta_{n + 1}^2} \right)^{1/2} (b\ci I \mp \delta_{n + 1} a\ci I  ) ; 
\end{align}
note that $a\ci{I_+}$, $b\ci{I_+}$ and $a\ci{I_-}$, $b\ci{I_-}$ are orthonormal pairs. 

Let us now find the relations between $\alpha_n$, $\beta_n$,  $\delta_{n+1}$, $\alpha_{n+1}$, and 
$\beta_{n+1}$, so that the martingale dynamics  \eqref{e: mart dyn} holds. We have 
%
\begin{eqnarray*}
 W \ci I & = & \frac12 \left( W\ci{I_-} +  W \ci{I_+} \right) \\
  & = & \frac12\cdot \frac{\alpha_{n + 1}}{1 + \delta_{n + 1}^2} (a\ci I -
  \delta_{n + 1} b\ci I) (a\ci I^{\ast} - \delta_{n + 1} b\ci I^{\ast})\\
  &  & + \frac12\cdot  \frac{\beta_{n + 1}}{1 + \delta_{n+1}^2} (\delta_{n + 1} a\ci I
  + b\ci I) (\delta_{n + 1} a\ci I^{\ast} + b\ci I^{\ast})\\
  &  & + \frac12\cdot  \frac{\alpha_{n + 1}}{1 + \delta_{n + 1}^2} (a\ci I +
  \delta_{n + 1} b\ci I) (a\ci I^{\ast} + \delta_{n + 1} b\ci I^{\ast})\\
  &  & + \frac12\cdot  \frac{\beta_{n + 1}}{1 + \delta_{n + 1}^2} (- \delta_{n
  + 1} a\ci I + b\ci I) (- \delta_{n + 1} a\ci I^{\ast} + b\ci I^{\ast})\\
  & = & \frac{\alpha_{n + 1} + \beta_{n + 1} \delta_{n + 1}^2}{1 + \delta_{n
  + 1}^2} a\ci I a\ci I^{\ast}\\
  &  & + \frac{\alpha_{n + 1} \delta_{n + 1}^2 + \beta_{n + 1}}{1 + \delta_{n
  + 1}^2} b\ci I b\ci I^{\ast} . \label{coefficient-recursion}
\end{eqnarray*}
Thus, the martingale dynamics \eqref{e: mart dyn} holds if and only if 
\begin{equation}
  \alpha_n = \frac{\alpha_{n + 1} + \beta_{n + 1} \delta_{n + 1}^2}{1 +
  \delta_{n + 1}^2} \qquad \text{and}\qquad \beta_n = \frac{\alpha_{n + 1} \delta_{n + 1}^2
  + \beta_{n + 1}}{1 + \delta_{n + 1}^2} \label{recursion-fact} .
\end{equation}
Note that it follows from relations \eqref{recursion-fact} that 
\begin{align}
\label{e: alpha+beta}
\alpha_n+\beta_n = \alpha_{n + 1} + \beta_{n+1}, 
\end{align}
as it should be, according to the martingale dynamics \eqref{e: mart dyn} and the linearity of the trace. 

Now, given $\alpha_n$, $\beta_n$, $\alpha_{n+1}$, $\beta_{n+1}$ satisfying \eqref{e: alpha+beta},  
we can easily find $\delta_{n+1}$ by solving equations \eqref{recursion-fact}, 
\begin{align}
\label{e: delta^2}
\delta_{n+1}^2 
= \frac{\alpha_{n+1} - \alpha_n}{\alpha_n - \beta_{n+1}}
= \frac{\beta_n - \beta_{n+1}}{\alpha_{n+1} - \beta_n} ;
\end{align}
the two expressions for $\delta_{n+1}^2$ here coincide because of \eqref{e: alpha+beta}. 

We now want the condition numbers $\alpha_n/\beta_n$ to increase exponentially, so let us take  
\begin{equation}
\label{eccentricity}
\frac{\beta_n}{\alpha_n} = \varepsilon^{2 n + 2}  , \qquad n=0, 1, 2, \ldots \ ,  
\end{equation}
where a small $\e>0$ is to be chosen later. For the sake of  convenience in the calculations, we want all 
$\delta_n$ to be small, and the extra $\e^2$ in \eqref{eccentricity} helps with that.

If we fix the sums in \eqref{e: alpha+beta}, 
$\alpha_n+\beta_n=1$, then we get from  \eqref{eccentricity} that 
%
%
\begin{equation}
\alpha_n = \frac{1}{1 + \varepsilon^{2 n + 2}} , \qquad \text{and}\qquad \beta_n =
\frac{\varepsilon^{2 n + 2}}{1 + \varepsilon^{2 n + 2}}
\label{coefficientsalphabetaofepsilon} \, .
\end{equation}
Substituting these expressions into \eqref{e: delta^2} (with $n$ replaced by 
$n-1$), we get 
\begin{align}
\label{e: delta_n^2}
\delta^2_n = \frac{\varepsilon^{2 n} (1 - \varepsilon^2)}{1 -
	\varepsilon^{4 n + 2}}, \qquad n=1, 2, 3, \ldots  
\end{align}

So, as we discussed above in Section \ref{s: overview no CET}, there exists a weight $W$ with its 
averages given by $\langle W \rangle\ci I = W\ci I$ for all $I\in\cD$.

\subsubsection{$A\ci I$ and the testing condition \eqref{e: Carl test 02}}

We define $A\ci I^{1/2}$ to be multiples of the projections onto $b\ci I$, i.e.\hspace{-0.05em}, 
\[ 
A^{1/2}\ci{I} = | I |^{1/2} r\ci I b\ci I b^{\ast}\ci I, 
\]
where for $I\in\cD^n$
\[ 
r\ci I = r_n = \frac{1}{\varepsilon^{n + 1}} . 
\]
First, we will prove that the testing condition \eqref{e: Carl test 02} holds, i.e.\hspace{-0.05 em}, that for every 
dyadic interval $K \in \mathcal{D}$ and every  vector $e\in\R^2$ we have
\begin{equation}
\label{CarlesonCondition}
  \sum_{I \in \mathcal{D} (K)} \lnm A\ci I^{1/2} \langle W e \rangle\ci I
  \rnm\ci{\mathbb{R}^2}^2  \leqslant C |K| \left\langle \langle W \rangle\ci K e, e   
  \right\rangle\ci {\R^2}.
\end{equation}
Let $K \in \mathcal{D}$ be such that $| K | = 2^{- n_0}$ for some $n_0
\geqslant 0$.  Let us denote the sum on the left hand side of 
\eqref{CarlesonCondition} by  $\Sigma_1$. Observe that
\[ 
A^{1/2}\ci{I} \langle W \rangle\ci I = | I |^{1/2} \beta\ci I r\ci I b\ci I b^{\ast}\ci I = 
\beta\ci I A\ci I^{1/2}, 
\]
where we recall that $\beta\ci I = \beta_n$ for $I\in\cD^n$. So we get
\begin{eqnarray*}
  \Sigma_1 & = & \sum_{I \in \mathcal{D} (K)} 
  \lnm A\ci I^{1/2} \langle W e \rangle\ci I \rnm\ci{\mathbb{R}^2}^2 
  = \sum_{I \in \mathcal{D} (K)}
   \lnm A\ci I^{1/2} \langle
  W \rangle\ci I e \rnm\ci{\mathbb{R}^2}^2\\
  & = & \sum_{I \in \mathcal{D} (K)} | I | (\beta\ci I r\ci I)^2 \lnm b\ci I b^{\ast}\ci I e 
  \rnm\ci{\mathbb{R}^2}^2 \\
  & = & \sum_{I \in \mathcal{D} (K)} | I | (\beta\ci I r\ci I)^2  \langle e, b\ci I
  \rangle\ci{\mathbb{R}^2} ^2 .
\end{eqnarray*}
Decomposing 
\begin{align*}
e= e_1 a\ci K + e_2 b\ci K, \qquad e_1,\,e_2\in\R, 
\end{align*}
we see that 
\begin{align*}
 \langle e, b\ci I \rangle\ci{\mathbb{R}^2} ^2 \leqslant \lnm e \rnm\ci{\R^2}^2 = e_1^2 + e_2^2. 
\end{align*}
Using the fact that for $I\in\cD^n$ we have $\beta\ci I = \beta_n \leqslant
\varepsilon^{2 n + 2} $ and $r\ci I = r_n = \frac{1}{\varepsilon^{n + 1}}$, we estimate
\begin{eqnarray*}
  \Sigma_1 & \leqslant & | e |\ci{\mathbb{R}^2}^2 \sum_{I \in \mathcal{D} (K)}
  | I | (\beta\ci I r\ci I)^2 
\leqslant | e |\ci{\mathbb{R}^2}^2 
\sum_{n = n_0}^{\infty}  2^{- n_0} \left(\varepsilon^{2 n + 2} \cdot
  \varepsilon^{- n - 1} \right)^2\\
  & = & \lnm e \rnm\ci{\mathbb{R}^2}^2 2^{- n_0} \varepsilon^{2 n_0 + 2}  \sum_{n =
  0}^{\infty} \varepsilon^{2 n} = 2^{- n_0} \frac{\varepsilon^{2 n_0 +
  2}}{1 - \varepsilon^2} \left(e_{1^{}}^2 + e_2^2 \right) .
\end{eqnarray*}
The right-hand side of (\ref{CarlesonCondition}) can be estimated from below
\begin{align*}
| K | \left\langle \langle W \rangle\ci K e, e \right\rangle\ci {\mathbb{R}^2} 
  & =  2^{- n_0} (\alpha_{n_0} e_1^2 + \beta_{n_0} e_2^2)   \\
  & \geqslant
  2^{- n_0} (1 - \varepsilon^2) (e_1^2 + \varepsilon^{2 n_0 + 2} e_{2^{}}^2) \\
  & \geqslant 2^{- n_0} (1 - \varepsilon^2)  \varepsilon^{2 n_0 + 2} (e_1^2 +  e_{2^{}}^2) .
\end{align*}
In the second inequality we used the fact that $\alpha_{n_0} \geqslant 1 -
\varepsilon^2$ and $\beta_{n_0} \geqslant (1 - \varepsilon^2) \varepsilon^{2
n_0 + 2}$, derived from the equations in
(\ref{coefficientsalphabetaofepsilon}). 

Comparing this estimate with the above upper bound for $\Sigma_1$, we see that for a constant $C$
such that $C (1 - \varepsilon^2)^2 \geqslant 1$ we have 
\[ 
\Sigma_1 \leqslant C 
 |K| \left\langle \langle W \rangle\ci K e, e   \right\rangle\ci {\R^2},  
\]
so the testing condition \eqref{CarlesonCondition} holds. 
\subsubsection{The blow-up}
\label{s: blow-up}
Now we want to show that there exist a vector $e\in \R^2$ and scalar functions
$\varphi\ci I$ supported on $I$ with $- 1 \leqslant \varphi\ci I \leqslant 1$
such that for $f=\1\ci{I_0} e$, 
\begin{equation}
\label{BlowUp}
\Sigma_2:=\sum_{I \in \mathcal{D}} \lnm A\ci I^{1/2} \langle \varphi\ci I W
  f \rangle\ci I \rnm\ci{\mathbb{R}^2}^2  = \infty .
\end{equation}
Let us chose $e = a = a\ci{I_0}$.   Note that $\| \1\ci{I_0}e\|\ci{L^2_W(\R^2)} ^2= \langle W\ci {I_0} a,a \rangle \ci {\mathbb{R}^2}=\frac1{1+\varepsilon^2}<\infty$. 

We will choose $\varphi\ci I = \1\ci {I_+}$ for all intervals $I$. Therefore for $f= \1\ci{I_0}a$ we
have
\begin{equation}
\label{Qplus}
  \langle \varphi\ci I W f \rangle\ci I = \frac{1}{| I |} \int\ci I \varphi\ci I (x) W
  (x) a \, \mathd x = \frac{1}{| I |} \int_{I_+} W (x) a \, \mathd x = \frac{| I_+
  |}{| I |} \langle W \rangle\ci {I_+} a = \frac12 \langle W \rangle\ci{I_+}
  a  .
\end{equation}
Hence we can expand the sum in \eqref{BlowUp} as
\begin{eqnarray}
\notag
  \sum_{I \in \mathcal{D}} \lnm A\ci I^{1/2} \langle \varphi\ci I W
  a \rangle\ci I \rnm\ci{\mathbb{R}^2}^2 & = &\frac{1}{4} \sum_{I \in \mathcal{D}}  r\ci I ^2  |
  I | \lnm  b\ci I b^\ast\ci I\langle W \rangle\ci{I_+} a \rnm\ci{\mathbb{R}^2}^2\\
\label{e: blow 03}
  & = & \frac{1}{4} \sum_{I \in \mathcal{D}}  r\ci I ^2  | I | \, \left| \langle a, \langle
  W \rangle \ci{I_+} b\ci I \rangle\ci{\mathbb{R}^2} \right|^2 .
\end{eqnarray}
For $I\in\cD^n$ denote 
\begin{align*}
\gamma\ci I = \gamma_{n+1}= \arctan \delta_{n+1},  
\end{align*}
so the relations \eqref{e: ab to+-} can be rewritten as
\begin{align}
\label{e: ab via gamma}
a\ci{I_\pm} = (\cos\gamma\ci I)a\ci I  \pm  (\sin \gamma\ci I) b\ci I, \qquad 
b\ci{I_\pm} = (\cos\gamma\ci I)b\ci I \mp (\sin \gamma\ci I) a\ci I.
\end{align}
Then, since $\langle W \rangle\ci{I_+} = \alpha\ci{I_+}  a\ci{I_+}a\ci{I_+}^* +\beta\ci{I_+} 
b\ci{I_+} b\ci{I_+}^*$, we can see from \eqref{e: ab via gamma} that 
\[ 
\langle W \rangle\ci{I_+} b\ci I = \alpha\ci{I_+} \sin \gamma\ci{I} a\ci{I_+} + \beta\ci{I_+} 
\cos \gamma\ci {I} b\ci{I_+} . 
\]
Then we can rewrite $\Sigma_2$ from \eqref{e: blow 03} as 
\begin{equation}
\label{expr}
\frac{1}{4} \sum_{I \in \mathcal{D}}  r\ci I ^2   | I | \left( \alpha\ci{I_+} (\sin 
\gamma\ci{I} )
\langle a, a\ci{I_+} \rangle\ci{\R^2} + \beta\ci{I_+} (\cos \gamma\ci{I} )  \langle a, 
b\ci{I_+}\rangle\ci{\R^2} \right)^2 
=  \frac{1}{4} \sum_{I \in \mathcal{D}}  r\ci I ^2   | I | \left(   D\ci I + F\ci I  
\right)^2 ,  
\end{equation}
where 
\begin{align*}
D\ci I := \alpha\ci{I_+} (\sin \gamma\ci{I} ) \langle a, a\ci{I_+} \rangle\ci{\R^2}, 
\qquad
F\ci I := \beta\ci{I_+} (\cos \gamma\ci{I} )  \langle a, b\ci{I_+}\rangle\ci{\R^2} . 
\end{align*}

We will show that if  $\varepsilon$ is sufficiently small,  then 
\begin{align}
\label{e: F<D} 
D\ci I  & \geqslant 2 |F\ci I|  
\intertext{and}
\label{e: D > c/r_I}
D\ci I &\geqslant (8 r\ci I)^{-1}.
\end{align}
Then we can ignore terms $F\ci I$ in \eqref{expr}, and estimate $\Sigma_2$ (with some $c>0$)
\begin{align}
\label{e: final blow-up}
\Sigma_2\geqslant c \sum_{I\in\cD} r\ci I^2 |I| r\ci I ^{-2} = c\sum_{n=0}^\infty \sum_{I\in\cD^n} 
|I| =
c \sum_{n=0}^\infty 1 =\infty , 
\end{align}
thus proving \eqref{Qplus} modulo estimates \eqref{e: F<D} and \eqref{e: D > c/r_I}. 

So, let us prove estimates \eqref{e: F<D} and \eqref{e: D > c/r_I}. 
Trivially $|F\ci I |\leqslant \beta\ci{I_+} $, and since for $I\in\cD^n$ we have 
 $\beta\ci{I_+} = \beta_{n + 1} \leqslant \varepsilon^{2 n + 4} $,  it follows that
\begin{align}
\label{e: F le}
| F\ci I | \leqslant \beta\ci{I_+} \leqslant \varepsilon^{2n + 4} .
\end{align}

Now let us look at $D\ci I $.  Since 
\begin{align*}
0<\gamma_n = \arctan \delta_n \leqslant \delta_n \leqslant \e^n
\end{align*}
(the last inequality follows immediately from \eqref{e: delta_n^2}), 
the angle $\gamma$ between $a$ and $a\ci{I_+}$ can be bounded as
\begin{align}
\label{gamma le}
0\leqslant \gamma \leqslant \sum_{n=1}^\infty \gamma_n \leqslant \sum_{n=1}^\infty \e^n = 
\frac{\e}{1-\e}\,.
\end{align}
Recalling that $\lnm a\rnm\ci{\R^2} = | a\ci{I_+} | \ci{\R^2} =1$, we therefore can see that  
$\langle a, a\ci{I_+}\rangle\ci{\R^2} \geqslant 1/2$ for sufficiently small $\e$. 
Also, we have for sufficiently small $\e$
\begin{align*}
\sin \gamma\ci{I} & = \sin \gamma_{n+1} = \delta_{n + 1} (1 + \delta^2_{n + 1})^{- 1/2} \geqslant
\frac{\varepsilon^{n + 1}}{2},
\\
\alpha\ci{I_+} &= \alpha_{n + 1} \geqslant 1/2.
\end{align*}
Combining the above three estimates, we get that 
\begin{align}
\label{e: D ge 01}
D\ci I \geqslant \e^{n+1}/8 = (8r\ci I)^{-1}, 
\end{align}
i.e.\hspace{-0.05em}, that the estimate \eqref{e: D > c/r_I} holds (for sufficiently small $\e$). Comparing the 
above bound \eqref{e: D ge 01} with \eqref{e: F le}, we can immediately see that 
\eqref{e: F<D} holds for $0<\e \leqslant 1/2$. \hfill \qed

\begin{remark}
\label{half-blow-up}
Notice that the left-hand side of (\ref{BlowUp}) blows up even if we do not sum over all $I\in\cD$ but only over right half intervals, such as when 
$\varphi\ci I =  \1\ci {I_+}$ if $I \in \mathcal{D}_+$ and $\varphi\ci I = 0$ otherwise. 
In that case the estimates \eqref{e: F<D},  \eqref{e: D > c/r_I} still hold, so in the same way as before 
we get for $f= \1\ci{I_0}a$
\begin{align}
\label{e: final blow-up 1/2}
\sum_{I \in \mathcal{D}_+ } \lnm A\ci I^{1/2} \langle    
\varphi\ci I W f \rangle\ci I \rnm^2_{\mathbb{R}^2} \geqslant 
c \sum_{I\in\cD_+} r\ci I^2 |I| r\ci I ^{-2} = c\sum_{n=0}^\infty \sum_{I\in\cD^n_+} 
|I| = c \sum_{n=0}^\infty 1/2 =\infty   .  
\end{align}

\end{remark}

\

\subsection{Failure of Convex Body Weighted Maximal Theorem: proof of Theorem \ref{noWcMAX} }

\subsubsection{Construction}\label{s: constr}

Let us take the sequence $A\ci I$ and $W$ from the above example. We set, using
the constant from inequality (\ref{CarlesonCondition})
\begin{align}
\label{e: wt A_I}
\wt{A}\ci I = C^{- 1} \langle W \rangle\ci I A\ci I \langle W \rangle\ci I,
\end{align}
where due to the Carleson property in inequality (\ref{CarlesonCondition}) we
have
\begin{equation}
  \frac{1}{| I |} \sum_{J \in \mathcal{D} (I)} \wt{A}\ci J =
  \frac{1}{C | I |} \sum_{J \in \mathcal{D} (I)} \langle W \rangle\ci J A\ci J
  \langle W \rangle\ci J \leqslant \langle W \rangle\ci I . \label{carltilde}
\end{equation}
Fix $n$. Consider the weight $W_n$, 
\begin{align*}
W_n := \sum_{I\in\cD^n} \langle W \rangle\ci I \1\ci I; 
\end{align*}
so $W_n$ is just the martingale defining $W$ at time $n$. 

For an
interval $I\in\mathcal{D}^{\leqslant n}_+ $ denote by $S\ci I$ the leftmost 
interval of size $2^{- n - 1}$ contained in $I$. That is, the interval of size $2^{- n - 1}$ reached from $I$ via sign tosses to the left. Denote by $\mathcal{S}^n$ the  collection of all 
such intervals, $\cS^n:= \{S\ci I: I\in\cD^{\leqslant n}_+ \}$. Observe that the intervals in $\mathcal{S}^n$ are pairwise disjoint. Indeed, since they are all of equal length they are either disjoint or identical. But any interval $J$ in $\mathcal{S}^n$ cannot arise from both $I,I'$ in 
$\mathcal{D}^{\leqslant n}_+ $ with, say, $ J \subsetneq I'\subsetneq I$. By construction, the path from $I$ to $J$ consists only of sign tosses to the left. Since $I' \in \mathcal{D}^{\leqslant n}_+ $, the last sign toss from the path from $I$ to $I'$ was to the right, so it cannot lie on the path from $I$ to $J$.

Now we define a family of weights $ W_{n,s}$,  $0 \leqslant s $, 
\[ 
 {W}_{n, s} := W_n + s \sum_{I \in \mathcal{D}^{\leqslant n}_+ }|S\ci I|^{-1}\1\ci{ S\ci I} \wt{A}\ci 
 I,
\]
where the matrices $\wt A\ci I$ are defined above by \eqref{e: wt A_I}. Note that the weights $W_{n, s}$ are measurable in $\mathcal{F}_{n + 1}$, meaning that they are constant on intervals $I\in 
\cD^{n+1}$.

\
\begin{figure}[h]\centering
\includegraphics{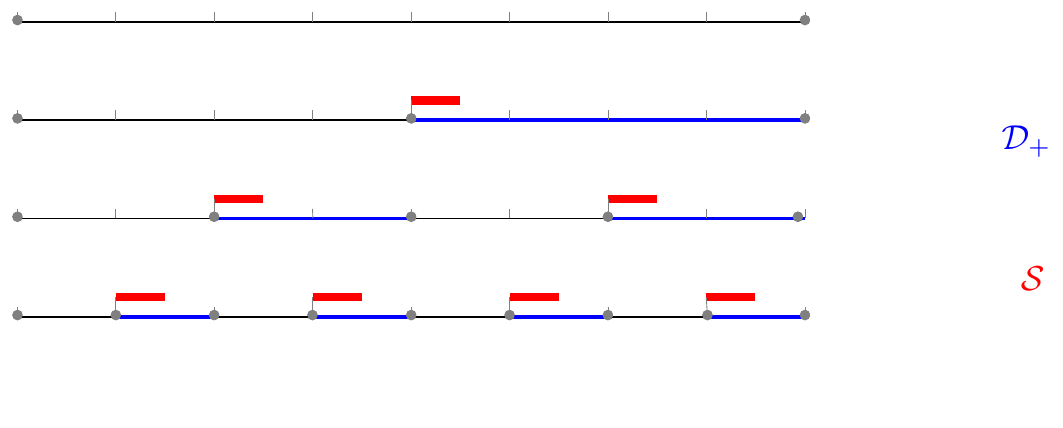}
\vspace{-0.8cm}
\caption{In this picture we present the construction of the intervals $S\ci I$ for $n=3$. The 
intervals $I\in\cD^{\leqslant n}_+$ are marked in blue, and the corresponding intervals $S\ci I$ 
are marked in red. }
\label{fig: S_I}
\end{figure}

To prove Theorem \ref{noWcMAX} we will show by contradiction  that for  the family of weights 
$W_{n,s}$ 
we \emph{do not have} the 
uniform estimate 
\begin{align}
\label{e: max bound}
\left\| M\ci{W_{n,s}}\ut{c} f  \right\|_{L^2(\R)}  \leqslant C \|f\|\ci{L^2\ci{W_{n,s}} (\R^2)} 
\qquad \forall f \in L^2\ci{W_{n,s}}(\R^2)
\end{align}
with $C$ not depending on $n$ and $s$. 
An elementary reasoning then gives us a weight $\wt W$ such that 
\begin{align*}
\Bigl\| M\ci{\wt W}\ut{c} \tilde f  \Bigr\|_{L^2(\R)} = \infty 
\end{align*}
for some function $\tilde f \in L^2_{\wt W} (\R^2)$. 

Indeed, let $s_k$, $n_k$ be such that for the weights $W_k:= W_{n_k, s_k}$ on $I_0$ there exist 
non-zero $f_k\in L^2\ci{W_k}(\R^2)$ supported on $I_0$ such that 
\begin{align}
\label{e: BlowUp 02}
\left\| M\ci{W_{k}}\ut{c} f  \right\|^2_{L^2}  > 4^k  \|f\|^2\ci{L^2\ci{W_k}(\R^2)}. 
\end{align}
Note, that this inequality is invariant under rescaling, meaning that it does not change if we 
multiply $W_k$ and $f_k$ by some non-zero constants. It is also easy to see that it does not change 
under an affine change of variables (applied to all objects simultaneously).  

To construct the weight $\wt W$ on $I_0$ let us represent $I_0$ as a union of disjoint intervals 
$I_k\in\cD$, $k\geqslant1$ (for example, take $I_k:= [2^{-k}, 2^{-k+1})$). Let $\wt W_k$ be 
the weight $W_k$ transplanted via an affine change of variables to the interval $I_k$ and normalized 
(by multiplying by a non-zero constant) in such a way that $\langle W_k \rangle\ci{I_k}\leqslant \bI$ (the identity matrix).  

Let also $\tilde f_k$ be the function $f_k$ transplanted by the same affine change as $W_k$ to the 
interval $I_k$ and normalized by $\|\tilde f_k\|^2\ci{\wt W_k}=2^{-k}$.  
 
Defining 
\[
\wt W (x) := \sum_{k=0}^\infty \1\ci{I_k}(x) \wt W_{k} (x) , \qquad 
\tilde f(x) := \sum_{k=0}^\infty \1\ci{I_k}(x) \tilde f_k(x),
\]
we immediately see that $\|\tilde f\|^2\ci{\wt W}=1$, and the estimates \eqref{e: BlowUp 02} imply 
that 
\[
\left\| \1\ci{I_k} M\ut{c}\ci{\wt W} \tilde f  \right\|^2_{L^2} \geqslant 
\left\| \1\ci{I_k} M\ut{c}\ci{\wt W_{k}} \tilde f_k  \right\|^2_{L^2}
> 4^k \|\tilde f_k\|^2\ci{\wt W_k} \geqslant 2^k,
\]
which  gives the desired blow-up.

\subsubsection{An a priori estimate}
Now, let us assume that we have the uniform estimate \eqref{e: max bound} just for one function 
$f=\1\ci{I_0}a $, where $a\in\R^2$ is the same as in inequality \eqref{e: final blow-up 1/2} in Remark \ref{half-blow-up}. 

Our goal is to arrive at a contradiction to \eqref{e: final blow-up 1/2}. From our assumption we first deduce some 
weaker estimate, from which using a trick from \cite{CuTr15} we will get the desired conclusion. 

To get to the final contradiction we need to estimate the weighted maximal function from below. We start 
with the trivial observation that for a weight $W$ and a function $f\in L^2\ci{W}(\R^2)$ and  any 
collection of disjoint measurable sets $S\ci I\subset I$ and measurable functions $\f\ci I:I\to 
[-1,1]$ 
parametrized by $I\in\cD$, the function $F=F[f,W, S\ci I, \f\ci I]$, 
\begin{align*}
F(x) := \sum_{I\in\cD} \lnm \1\ci{S_I}(x) W(x)^{1/2} \langle W \rangle\ci{I}^{-1} 
\langle \f\ci I W f \rangle\ci{I} \rnm\ci{\R^2}
\end{align*}
is pointwise estimated as $F(x) \leqslant M\ut c\ci W f (x)$. Indeed, if $x\notin \cup_{I\in\cD} S \ci I$, then $F(x)=0$. 
Otherwise the sum defining $F(x)$ consists of exactly one term
\begin{align}
\label{e: F(x)}
F(x) = \lnm W(x)^{1/2} \langle W \rangle\ci{I}^{-1} \langle \f\ci I W  f \rangle\ci{I} \rnm\ci{\R^2}
\end{align}
where $I=I(x)$ is the unique $I\in\cD$ such that $x\in S\ci I$; uniqueness of $I$ follows from the 
disjointness of $S\ci I$. To get the 
maximal function $M\ut c\ci Wf(x)$, we need to take the supremum of the right hand side of 
\eqref{e: F(x)} over all $I\in\cD$, $I\ni x$ and over all $\f\ci I$, so no particular choice of $I$ 
and $\f\ci I$ can give us more than the maximal function. 

It is easy to compute the norm of $F$ and rewrite the trivial inequality $\|F\|\ci{L^2(\R)}^2 \leqslant 
\|M\ut c\ci W f\|\ci{L^2(\R)}^2$ as 
\begin{align}
\label{e: LinMaxFn 01}
\sum_{I\in\cD} |S\ci I| \lnm  \langle W\rangle\ci{S_{{}_I}}^{1/2}   \langle W\rangle\ci I^{-1} 
\langle 
\f\ci I W 
f\rangle\ci I \rnm\ci{\R^2}^2 \leqslant \|M\ut c\ci W f \|\ci{L^2(\R)}^2 . 
\end{align}
The above inequality \eqref{e: LinMaxFn 01} holds for any collection of disjoint measurable 
sets $S\ci I\subset I$ and functions $\f\ci I: I\to[-1,1]$. 

We remark that if we take in the left hand side 
of \eqref{e: LinMaxFn 01} the supremum over all such collections, we will get exactly $ \|M\ut c\ci 
W f \|\ci{L^2(\R)}^2 $. This is well known to the experts as the \emph{linearization} for the 
maximal function, that reduces without loss of generality the estimates of a nonlinear operator 
(the maximal function) to the estimates of the (linear) embedding operator. Since for the current paper we only need the trivial estimate \eqref{e: LinMaxFn 01}, we will leave 
the proof of the details of this \emph{linearization} as an exercise for a curious reader.  

We specify the estimate \eqref{e: LinMaxFn 01} to the case of $W=W_{n,s}$ constructed above and $f= 
 \1\ci {I_0} a$,  with $S\ci I$ defined in Section \ref{s: constr}
and $\f\ci I=\1\ci{I_+}$ for all $I \in \cD^{\leqslant n}_+$ and $\varphi \ci I =0$ else, similarly as in Section \ref{s: blow-up}. By noticing 
that for $I\in\cD_+^{\leqslant n}$ we have 
\begin{align*}
W_{n,s} (S\ci I) =  W(S\ci I) + s \wt A\ci I \geqslant s \wt A\ci I, 
\end{align*}
we get from the estimate \eqref{e: LinMaxFn 01} that 
\begin{align}
\label{e: max sum 01}
\sum_{I\in\cD^{\leqslant n}_+}  s \lnm \wt A\ci I ^{1/2} \langle W_{n,s} \rangle\ci{I}^{-1} \langle 
\f\ci I W_{n,s} f \rangle\ci{I} \rnm\ci{\R^2}^2 
\leqslant 
\Bigl\|  M\ut c \ci{W_{n,s}} f  \Bigr\|\ci{L^2(\R)}^2 . 
\end{align}

But for $f=\1\ci{I_0}a$, $a\in\R^2$, we have 
\begin{align*}
\|f\|\ci{L^2_{W_{n,s}}}^2 =\left\langle \langle W \rangle\ci{I_0} a, a \right\rangle\ci{\R^2} +
s \sum_{I\in\cD_+} \left\langle  \wt A\ci I a, a \right\rangle\ci{\R^2}   
\leqslant (1+s) \left\langle \langle W \rangle\ci{I_0} a, a \right\rangle\ci{\R^2}  ,  
\end{align*}
so we get from \eqref{e: max sum 01} that  
\begin{align}
\label{e: max sum 02}
s \sum_{I\in\cD^{\leqslant n}_+} \lnm \wt A\ci I ^{1/2} \langle W_{n,s} \rangle\ci{I}^{-1} \langle 
\f\ci I W_{n,s} f  \rangle\ci{I} \rnm\ci{\R^2}^2 
\leqslant 
(1+s) C \left\langle \langle W \rangle\ci{I_0} a, a 
\right\rangle\ci{\R^2} .
\end{align} which is our a priori estimate.

\subsubsection{From estimate (\ref{e: max sum 02})  to a contradiction}

For $I\in \cD^{\leqslant n}_+$ define 
\begin{align*}
R\ci{n, I} (s):=  \lnm \wt A\ci I ^{1/2} \langle W_{n,s} \rangle\ci{I}^{-1} \langle 
\f\ci I W_{n,s} f  \rangle\ci{I} \rnm\ci{\R^2}^2 , 
\end{align*}
so  \eqref{e: max sum 02} can be rewritten as 
\begin{align}
\label{e: sum R_n}
s \sum_{I\in\cD^{\leqslant n}_+} R\ci{n, I} (s)  \leqslant   (1+s) C \left\langle \langle W 
\rangle\ci{I_0} 
a, a \right\rangle\ci{\R^2}   .
\end{align}

By the cofactor inversion formula, the entries of the matrix $\left\langle W_{n, 
s}\right\rangle_{I}^{-1}$ are rational functions of $s$ of the form $p\ci{n, I}(s) / Q\ci{n, I}(s)$ 
where $p\ci{n, I}(s)$ is affine in $s$ and $Q_{n, I}(s)=\det\left(\left\langle W_{n, 
s}\right\rangle_{I}\right)$ has degree $2$.

The components of the vector $\left\langle W_{n, s}\right\rangle\ci{I_{+}}\!\!\! a$ are polynomials 
of 
degree 1, so we can write
\begin{align*}
R_{n, I}(s)=\frac{P_{n, I}}{Q_{n, I}^{2}}(s), \quad Q_{n, 
I}(s)=\det\left(\left\langle W_{n, s}\right\rangle_{I}\right)
\end{align*}
and $P\ci{n, I}$ is a polynomial of degree at most $4$.

Recall that for $I\in\cD^{\leqslant n}_+$
\begin{align*}
\left\langle W_{n, s}\right\rangle\ci{I}=\langle W\rangle\ci{I}+s \frac{1}{|I|} \sum_{J \in 
\mathcal{D}_{+}^{\leqslant n}(I)} \tilde{A}_{J} \leqslant(1+s)\langle W\rangle\ci{I}
\end{align*}
and thus for all $s\geqslant 0$,
\begin{align*}
Q\ci{n, I}(s) \leqslant(1+s)^{2} \operatorname{det}\left(\langle W\rangle\ci{I}\right) = (1+s)^2 
Q\ci{n,I}(0).
\end{align*}
Therefore,  the estimate \eqref{e: sum R_n} implies that 
\begin{align}
\label{e: est p}
\sum_{I\in\cD^{\leqslant n}_+} \frac{P\ci{n, I} (s)}{Q\ci{n, I} (0)^2} \leqslant   \frac{(1+s)^5}{s} C 
\left\langle \langle W 
\rangle\ci{I_0} 
a,a \right\rangle\ci{\R^2} 
\end{align}

We need the following Lemma, see \cite[Lemma 2.2]{CuTr15}:

\begin{lemma}
	\label{lemma-polynomial}If $p$ is a polynomial such that 
	\[
	| p (s) |
	\leqslant \frac{(1 + s)^N}{s} \qquad \forall s > 0,
	\]
	then $| p (0) | \leqslant e^2N^2$.
\end{lemma}

We apply this lemma to the polynomial $p$, 
\begin{align*}
p(s) := \left( C \left\langle \langle W \rangle\ci{I_0} a,a \right\rangle\ci{\R^2} \right)^{-1} 
\sum_{I\in\cD^{\leqslant n}_+} \frac{P\ci{n,I}(s)}{Q\ci{n,I}(0)^2}. 
\end{align*}
Note that $P_{n,I}(s)$ are non-negative for $s \geqslant 0$.
Estimate \eqref{e: est p} means that $p$ satisfies the assumption of Lemma \ref{lemma-polynomial}  
with $N=5$, therefore
\begin{align}
\label{e: sum R_n (0) 02}
\sum_{I\in\cD^{\leqslant n}_+} R\ci{n, I}(0) & \leqslant C e^2 5^2 \left\langle \langle W 
\rangle\ci{I_0} a,a \right\rangle\ci{\R^2} .
\end{align}

Recall that for $I \in\cD^{\leqslant n}_+$ 
\begin{align*}
R_{n,I} (0)= \lnm \wt A\ci I ^{1/2} \langle W_{n} \rangle\ci{I}^{-1} \langle 
\f\ci I W_{n} f  \rangle\ci{I} \rnm\ci{\R^2}^2, 
\end{align*}
where $f=\1\ci{I_0} a$ and $\f\ci I = \1\ci{I_+}$. Noticing that for $I \in\cD^{\leqslant n-1}_+$ 
\begin{align*}
\langle W_n \rangle\ci I = \langle W \rangle\ci I, \qquad 
\langle \f\ci I W_n f\rangle\ci I =\langle \f\ci I W f\rangle\ci I , 
\end{align*}
we can deduce from \eqref{e: sum R_n (0) 02} that
\begin{align*}
\sum_{I\in\cD^{\leqslant n-1}_+} 
\lnm \wt A\ci I ^{1/2} \langle W \rangle\ci{I}^{-1} \langle \f\ci I W f  \rangle\ci{I} 
\rnm\ci{\R^2}^2
\leqslant C e^2 5^2 \left\langle \langle W \rangle\ci{I_0} a, a \right\rangle\ci{\R^2} .
\end{align*}
Letting $n\to\infty$ we get that 
\begin{align}
\label{e: blow-up bound 01}
\sum_{I\in\cD_+} 
\lnm \wt A\ci I ^{1/2} \langle W \rangle\ci{I}^{-1} \langle \f\ci I W f  \rangle\ci{I} 
\rnm\ci{\R^2}^2
\leqslant C e^2 5^2 \left\langle \langle W \rangle\ci{I_0} a, a \right\rangle\ci{\R^2} .
\end{align}
Recall that by \eqref{e: wt A_I} we have $\wt A\ci I = C^{-1} \langle W \rangle \ci I A\ci I 
\langle W \rangle \ci I$, so 
we rewrite the above estimate \eqref{e: blow-up bound 01} as 
\begin{align*}
\sum_{I\in\cD_+} 
\lnm  A\ci I ^{1/2}  \langle \f\ci I W f  \rangle\ci{I} 
\rnm\ci{\R^2}^2
\leqslant C e^2 5^2 \left\langle \langle W \rangle\ci{I_0} a, a \right\rangle\ci{\R^2} < \infty, 
\end{align*}
which contradicts the blow-up estimate \eqref{e: final blow-up 1/2} obtained before. 
\hfill \qed

\end{document}